\documentclass[11pt,a4paper]{article}
\usepackage[T1]{fontenc}
\usepackage{lmodern,microtype,needspace}
\usepackage[margin=28mm]{geometry}
\usepackage{amsmath,amssymb,amsthm,mathtools,booktabs,array}
\usepackage[hidelinks]{hyperref}
\hypersetup{pdftitle={Uniform cubic bounds for irrational values of the dilogarithm},pdfauthor={Hu Tan and Ying Zhang}}
\newcommand{\Li}{\operatorname{Li}_2}
\newcommand{\lcm}{\operatorname{lcm}}
\newcommand{\Z}{\mathbb Z}
\newcommand{\Q}{\mathbb Q}
\newcommand{\R}{\mathbb R}

\newcommand{\dd}{\,\mathrm d}
\newcommand{\pos}[1]{\left(#1\right)_+}
\newtheorem{theorem}{Theorem}[section]
\newtheorem{proposition}[theorem]{Proposition}
\newtheorem{lemma}[theorem]{Lemma}

\theoremstyle{remark}
\numberwithin{equation}{section}
\title{Uniform cubic bounds for irrational values of the dilogarithm}
\author{Hu Tan\thanks{Academy of Mathematics and Systems Science, Chinese Academy of Sciences, 55 Zhongguancun East Road, Beijing 100190, China. Email: \texttt{tanhu2020@amss.ac.cn}.}
\and Ying Zhang\thanks{School of Mathematical Sciences, Soochow University, 1 Shizi Street, Suzhou 215006, China. Email: \texttt{yzhang@suda.edu.cn}. Corresponding author.}}
\date{}
\begin{document}
\maketitle
\begin{abstract}
We study explicit denominator conditions ensuring the irrationality of the dilogarithm at positive and negative rational points. We prove that $\Li(p/q)$ is irrational whenever $q>(4e^4/27)p^3$, and obtain an explicit, slightly smaller threshold for $\Li(-p/q)$. The positive argument specializes the Rhin--Viola integral family and uses a uniform estimate for a restricted prime divisor. The negative argument combines the same arithmetic construction with three adjacent integral moments and an explicit maximization. We also give finite-numerator refinements by interpolating certified partition bounds and by introducing a comparison measure with a quadratic tail. In particular, $\Li(5/q)$ is irrational for every integer $q\ge418$. The proof separates the established integral and determinant constructions from the parameter estimates and finite certificates used here.
\end{abstract}

\section{Introduction}
The aim of this paper is to give explicit uniform denominator bounds for irrational values of the dilogarithm, and to refine them when the numerator is small.
For $|x|\le1$, we use $\Li(x)=\sum_{m\ge1}x^m/m^2$. The classical values at $1$ and $-1$ are $\pi^2/6$ and $-\pi^2/12$. Away from these special points, a convergent defining series does not itself provide an effective irrationality argument: one needs integer linear forms whose arithmetic denominators grow more slowly than their analytic decay.

Integral constructions have been central to this problem since Beukers's proof for $\zeta(2)$ and $\zeta(3)$ \cite{Beukers}. Hata \cite{Hata} developed rational approximations to the dilogarithm and proved, in particular, irrationality at $1/s$ for integers $s\ge7$ and $s\le-5$. Rhin and Viola \cite{RV05} introduced a five-parameter family of double integrals together with a permutation group that supplies additional common divisors of their coefficients. Their results give computable thresholds for arbitrary fixed positive numerators; the first two include $\Li(1/s)\notin\Q$ for $s\ge6$ and $\Li(2/s)\notin\Q$ for $s\ge51$. Rhin and Viola subsequently extended the construction to negative arguments and to linear independence of $1$, $\operatorname{Li}_1$, and $\Li$ \cite{RV19}. We use their polynomial identities and divisibility results as established arithmetic input.

Recent work on individual small rational points provides further useful input. Prellberg \cite{Prellberg} proves irrationality at $-1/4$, $1/5$ and $-1/3$. His negative-argument proof uses three bounded parameter shifts and a four-term moment recurrence to resolve nonvanishing. We apply that strategy to the symmetric Legendre family, with a variable negative argument, and obtain the explicit uniform threshold in Theorem~\ref{thm:negative}.

The fixed-numerator arguments follow a separate recent line of work. Fauzan's proof of the irrationality of $\zeta(5)$ \cite{Fauzan} uses normalized Hankel determinants, moment representations and logarithmic-energy estimates. Qian Tang adapted this approach to the dilogarithm, including the value $\Li(1/2)$ \cite{Tang}, and Si-Qi Liu extended the method to further rational arguments in \texttt{Li2Rational} \cite{Liu}. This mathematical lineage is recorded in the latter project's \texttt{README} and \texttt{NOTICE}. Our small-numerator refinements, including the bound $q\ge418$ for $\Li(5/q)$, build on the integer polynomials, arithmetic estimates, nonvanishing theorem and finite-circle comparison inequality in the pinned version of \cite{Liu}. We retain its normalization and supply the interpolation estimates and quadratic comparison measure used below.

Our uniform statement has a simple cubic scale. Put $C_0=4e^4/27=8.088614819725\ldots$. All numerators and denominators in the following statements are positive integers; coprimality is not required.

\begin{theorem}\label{thm:positive}
If $q>C_0p^3$, then $\Li(p/q)$ is irrational.
\end{theorem}

\begin{theorem}\label{thm:negative}
Define
\begin{equation}\label{eq:negative-threshold}
B_-(p)=\frac{2e^4p^3}{27}\left[1+\left(1+\frac{3}{e^4p^2}\right)^{3/2}\right]-\frac{2p}{3}.
\end{equation}
If $q>B_-(p)$, then $\Li(-p/q)$ is irrational. In particular, each of the conditions $q>C_0p^3-p/3+1/(4e^4p)$ and $q>C_0p^3$ is sufficient.
\end{theorem}

The negative threshold satisfies $B_-(p)=C_0p^3-p/3+1/(4e^4p)-1/(8e^8p^3)+O(p^{-5})$. This expansion describes the size of the improvement; the sufficient conditions follow from a global inequality, proved in Section~\ref{sec:negative}, rather than from the expansion alone.

\begin{theorem}\label{thm:finite}
For each row of Tables~\ref{tab:thresholds} and~\ref{tab:additional}, $\Li(p/q)$ is irrational for every integer $q\ge Q_+$, and $\Li(-p/q)$ is irrational for every integer $q\ge Q_-$.
\end{theorem}
\begin{table}[ht]
\centering
\begin{tabular}{rrrrrrrr}\toprule
$p$&1&2&3&4&5&6&7\\\midrule
$Q_+$&1&12&54&168&418&1038&2233\\
$Q_-$&1&16&60&178&434&1058&2256\\\bottomrule
\end{tabular}
\caption{Uniform sufficient thresholds for fixed numerators.}\label{tab:thresholds}
\end{table}

These thresholds apply to entire denominator ranges. Exact cancellation also gives additional isolated cases below them and includes the first denominator in several rows of the table.

\subsection*{Proof strategy and contributions}
The proof has three parts. For positive arguments we fix the parameter family $(2K,K-1,K,K,K)$ in the Rhin--Viola construction. A restricted first-period prime divisor suffices: its logarithmic gain is at least $(\tfrac32\log K-1)/K$ for $K\ge32$. Taking $K=2p$ proves the large-numerator range, while four rational certificates cover $3\le p\le79$. This gives one explicit constant for all numerators.

For negative arguments, we use three adjacent moments of a Legendre-type function. A four-term recurrence and uniqueness of polynomial moments ensure that these three moments cannot all vanish. Their common exponential majorant reduces the argument to maximizing $Ax(1-x)^2-x$ on $[0,1]$, where $A=e^4p^2$. The resulting maximum is the expression in \eqref{eq:negative-threshold} divided by $p$.

For small numerators, positive coefficients make the logarithm of the absolute partition sum convex in the logarithm of its argument. This allows interpolation between existing certified profiles. A new comparison density, constant on an initial interval and quadratic on its remaining support, gives the rate needed for the last three rows of Table~\ref{tab:thresholds}. We provide the approximation argument connecting this continuous density to the existing finite-circle theorem, with the comparison family fixed before the determinant order tends to infinity.

\section{Arithmetic input and irrationality criteria}\label{sec:input}
This section states the established arithmetic constructions in the normalizations used below and records the elementary final contradiction.
Write $d_m=\lcm(1,\ldots,m)$, with $d_0=1$. The prime number theorem gives $\log d_m=m+o(m)$; see, for example, \cite[pp.~111--114]{Davenport}.

\begin{lemma}\label{lem:linearforms}
Let $\xi\in\R$. Suppose that, for each sufficiently large $n$, finitely many numbers $L_{n,r}$ belong to $\Z+\Z\xi$, their maximum absolute value tends to zero, and at least one is nonzero. Then $\xi$ is irrational.
\end{lemma}
\begin{proof}
If $\xi=a/b$ with $a\in\Z$ and $b\in\Z_{>0}$, then each $bL_{n,r}$ is an integer. A nonzero one has absolute value at least $1$, contradicting their common convergence to zero.
\end{proof}

For nonnegative integer parameters satisfying the usual four nonnegativity conditions, let $I_z(h,j,k,l,m)$ denote the Rhin--Viola linear combination of a real double integral and its mixed residue integral; we follow \cite[Section~2]{RV05} for $z>1$ and \cite[Section~2]{RV19} for $z<0$. Their polynomial identity has the form
\begin{equation}\label{eq:rv-input}
d_Hd_K z^\alpha(z-1)^\beta I_z=P(z)-Q(z)\Li(1/z),\qquad P,Q\in\Z[z],
\end{equation}
with explicit degree bounds. We use the notation for $H,K,\alpha,\beta,\delta$ from those sources only when specializing this identity; the letter $K$ in Section~\ref{sec:positive} is our freely chosen family parameter. The permutation action gives further divisibility by products of primes detected by floor-function inequalities \cite[Lemma~4.1 and (4.4)]{RV05}.

The fixed-numerator argument uses a different established construction. Denote its positive absolute partition sum by $\mathcal Z_n(t)$, where $0<t\le1/2$, and put
\begin{equation}\label{eq:arithmetic-rate}
R_\sigma(p,q)=\frac83\log q-\frac94\log(q-\sigma p)-6\log p+\frac{72}{25},\qquad \sigma\in\{-1,1\}.
\end{equation}
The determinant and arithmetic results of \cite{Liu} yield the following useful formulation.

\begin{proposition}[Partition criterion]\label{prop:partition}
Suppose $0<t\le1/2$ and $\mathcal Z_n(t)\le\exp(Bn^2)$ for all sufficiently large $n$. If $p/q\le t$ and
\begin{equation}\label{eq:partition-gap}
B+\log\frac{p}{qt}<R_\sigma(p,q),
\end{equation}
then $\Li(\sigma p/q)$ is irrational.
\end{proposition}
\begin{proof}
The node-scaling inequality is $\mathcal Z_n(p/q)\le(p/(qt))^{n(n+1)}\mathcal Z_n(t)$. The integer polynomials in the cited construction have degree at most $2n$, and their arithmetic growth is governed by \eqref{eq:arithmetic-rate}. The strict gap in \eqref{eq:partition-gap} therefore gives integer polynomials $P_n$ with $|P_n(\Li(\sigma p/q))|\le e^{-cn^2}$ for some $c>0$, and infinitely many of these values are nonzero. If $\Li(\sigma p/q)=u/v$, then $v^{2n}P_n(u/v)$ is an integer and its absolute value tends to zero along that infinite subsequence, a contradiction. For unreduced $p/q$, cancelling a common factor $h$ increases the arithmetic rate by $(67/12)\log h$, so the stated raw rate remains sufficient.
\end{proof}

The polynomial construction, arithmetic growth estimate and nonvanishing theorem in this proposition are inputs from \cite{Liu}. The interpolation and comparison estimates supplying $B$ are proved below.

\section{The positive cubic bound}\label{sec:positive}
The aim of this section is to prove Theorem~\ref{thm:positive} by choosing one family of Rhin--Viola parameters and estimating its prime gain uniformly.
Fix $p,q$ and put $z=q/p>1$. Choose an integer $K\ge2$, write $\varepsilon=1/K$, and set $N=Kn$. The parameter $K$ is fixed while odd positive integers $n$ tend to infinity. Specialize the five parameters to $(2N,N-n,N,N,N)$.

\subsection{Integral forms and their sign}
Removing the factor $z^{-2N}$ from the normalization in \cite{RV05}, the real integral is
\begin{equation}\label{eq:J0}
J_{0,n}=\int_0^1\!\int_0^1
\frac{x^{N-n}(1-x)^{2N}y^N(1-y)^N}{[x+(z-x)y]^{N-n+1}}\dd y\dd x.
\end{equation}
Let $J_{1,n}$ be obtained by replacing the inner integration by the normalized, positively oriented residue around its unique pole, and put $F_n=J_{0,n}-(\log z)J_{1,n}$. If $r=N-n$, expansion at the pole gives
\begin{equation}\label{eq:residue}
J_{1,n}=(-1)^n\int_0^1\frac{x^N(1-x)^{2N}}{(z-x)^{N+1}}
\sum_{i=0}^N\binom Ni\binom{N+i}{r}\left(\frac{x}{z-x}\right)^i\dd x.
\end{equation}
Every factor inside the integral is nonnegative, and the integrand is positive on $(0,1)$. Hence $J_{1,n}<0$ for odd $n$, and $F_n>0$.

For these parameters, the invariants in the divisibility theorem are $M=N_{\mathrm{RV}}=\alpha=\delta=2K$ and $\beta=0$ before multiplication by $n$. Let $\Omega_K\subset(0,1)$ consist of the points satisfying at least one of
\begin{align}
\lfloor(2K-1)\omega\rfloor+\lfloor(K+1)\omega\rfloor
&<\lfloor2K\omega\rfloor+\lfloor K\omega\rfloor,\label{eq:floorA}\\
\lfloor(K-1)\omega\rfloor+\lfloor(K+1)\omega\rfloor
&<2\lfloor K\omega\rfloor.\label{eq:floorB}
\end{align}
Let $\Delta_n$ be the product of the primes $\ell>\sqrt{2Kn}$ with $n/\ell\in\Omega_K$. The divisibility result in \cite{RV05} and clearing the denominator of $z=q/p$ give
\begin{equation}\label{eq:positive-integer}
L_n=\frac{d_{2N}^2}{\Delta_n}p^{2N}F_n\in\Z+\Z\Li(p/q).
\end{equation}
Only the indicated first-period primes are used. Since $\Omega_K$ is a finite union of rational intervals bounded away from zero, the prime number theorem applied to each interval gives
\begin{equation}\label{eq:gain}
\lim_{n\to\infty}\frac{\log\Delta_n}{Kn}=G_K,
\qquad G_K=\frac1K\int_{\Omega_K}\frac{\dd\omega}{\omega^2}.
\end{equation}
Changing any of the finitely many interval endpoints does not change this limit.

\subsection{A uniform lower bound for the prime gain}
\begin{lemma}\label{lem:gain}
For every integer $K\ge32$, $KG_K>\tfrac32\log K-1$.
\end{lemma}
\begin{proof}
Put $m=\lceil K/4\rceil$ and change variables to $u=K\omega$. Direct substitution into \eqref{eq:floorA}--\eqref{eq:floorB} shows that the following disjoint half-open intervals are admissible:
\begin{equation}\label{eq:gain-intervals}
\left[k,\frac{k}{1-\varepsilon}\right)\ (1\le k\le m),
\qquad
\left[k+\frac12,\frac{k+1/2}{1-\varepsilon/2}\right)\ (0\le k<m).
\end{equation}
Indeed, throughout these intervals $\varepsilon u<1/3$. In the first family $u$ has integer part $k$ while $(1-\varepsilon)u<k$ and $(1+\varepsilon)u<k+1$, which proves \eqref{eq:floorB}. In the second family $2u$ has integer part $2k+1$, whereas $(2-\varepsilon)u<2k+1$ and $(1+\varepsilon)u<k+1$, which proves \eqref{eq:floorA}. The same inequalities ensure that an interval in the first family ends before the next half-integer and an interval in the second family ends before the next integer.

Writing $H_m=\sum_{j=1}^m1/j$, integration over \eqref{eq:gain-intervals} yields $KG_K\ge\tfrac12H_m+H_{2m}$. The elementary estimate $H_j\ge\log(j+1)+1-\log2$ gives
\begin{align*}
KG_K&\ge\tfrac12\log(m+1)+\log(2m+1)+\tfrac32(1-\log2)\\
&\ge\tfrac32\log K+\tfrac32-\tfrac72\log2
>\tfrac32\log K-1.
\end{align*}
Here $m+1\ge K/4$ and $2m+1\ge K/2$ for $K\ge32$; the last inequality follows, for example, from $\log2<5/7$.
\end{proof}

\subsection{Analytic rates}
Set $a_\varepsilon=4(1-\varepsilon)^{1-\varepsilon}/(3-\varepsilon)^{3-\varepsilon}$, $b_\varepsilon=\varepsilon^\varepsilon/(1+\varepsilon)^{1+\varepsilon}$ and $h(\varepsilon)=-\varepsilon\log\varepsilon-(1-\varepsilon)\log(1-\varepsilon)$.
Since $x+(z-x)y\ge zy$, the exponential part of \eqref{eq:J0} is bounded by
$z^{-1+\varepsilon}x^{1-\varepsilon}(1-x)^2y^\varepsilon(1-y)$.
Its two one-variable maxima are $a_\varepsilon$ and $b_\varepsilon$. The remaining weight $[x+(z-x)y]^{-1}$ is integrable on the square: near the only possible singular corner its denominator is comparable with $x+y$. Consequently
\begin{equation}\label{eq:J0-rate}
\limsup_{n\to\infty}J_{0,n}^{1/N}\le z^{-1+\varepsilon}a_\varepsilon b_\varepsilon.
\end{equation}

For the residue integral, $\binom{N+i}{N-n}/\binom N{N-n}=\prod_{j=1}^i(N+j)/(n+j)\le K^i$. Applying the binomial theorem in \eqref{eq:residue}, using $\max_{[0,1]}x(1-x)^2=4/27$, and using the finite entropy bound $\binom N{N-n}\le\exp(Nh(\varepsilon))$ gives
\begin{equation}\label{eq:J1-rate}
\limsup_{n\to\infty}|J_{1,n}|^{1/N}
\le\frac4{27(z-1)}\exp(h(\varepsilon))\left(1+\frac K{z-1}\right).
\end{equation}
Together with \eqref{eq:gain} and $\log d_{2Kn}^2/(Kn)\to4$, these estimates show that $L_n\to0$ provided both
\begin{align}
E_0&=4-G_K+2\log p-(1-\varepsilon)\log z+\log a_\varepsilon+\log b_\varepsilon,\label{eq:E0}\\
E_1&=4-G_K+2\log p+\log(4/27)-\log(z-1)
+h(\varepsilon)+\log\left(1+\frac K{z-1}\right)\label{eq:E1}
\end{align}
are negative. For fixed $p,K$, both expressions decrease as $z$ increases. It therefore suffices to check $z=C_0p^2$.

\subsection{Large and small numerators}
For $p\ge80$, take $K=2p$. The elementary logarithmic estimates
$\log a_\varepsilon\le\log(4/27)+\varepsilon\log4$,
$\log b_\varepsilon\le-\varepsilon\log K-\varepsilon$, and
$h(\varepsilon)\le\varepsilon\log K+\varepsilon$ apply. The second follows from $(1+\varepsilon)\log(1+\varepsilon)\ge\varepsilon$, and the third from $-(1-\varepsilon)\log(1-\varepsilon)\le\varepsilon$. For the first, differentiate the logarithm of $a_\varepsilon$: its derivative is $\log((3-\varepsilon)/(1-\varepsilon))\le\log4$ when $0\le\varepsilon\le1/3$.
Lemma~\ref{lem:gain} now gives
\begin{align}
E_0&\le\frac{-\tfrac12\log p+\log(C_0/\sqrt2)}{2p}<0,\label{eq:tail0}\\
E_1&\le\frac{-\tfrac12\log(2p)+2+\dfrac{2p(1+2p)}{8p^2-1}}{2p}
<\frac{-\tfrac12\log160+2.51}{2p}<0.\label{eq:tail1}
\end{align}
The first strict inequality uses $C_0<9$, and the second uses $C_0>8$, $\log(1+s)\le s$, and $p\ge80$. All constants have substantial rational margins; the accompanying exact checks use $8.08<C_0<8.09$.

For $3\le p\le79$, enumerate the rational breakpoints of the floor functions in \eqref{eq:floorA}--\eqref{eq:floorB}, and integrate $\omega^{-2}$ on the accepted intervals. Table~\ref{tab:positive-certificates} gives the resulting lower bounds and negative upper bounds at $z=C_0p^2$.
\begin{table}[ht]
\centering\small
\begin{tabular}{rrrrr}\toprule
$p$&$K$&lower bound for $G_K$&upper bound for $E_0$&upper bound for $E_1$\\\midrule
$3$--$4$&6&0.547222&$-0.0220$&$-0.0024$\\
$5$--$12$&14&0.331074&$-0.0089$&$-0.0014$\\
$13$--$58$&47&0.137880&$-0.00050$&$-0.00030$\\
$59$--$79$&114&0.068567&$-0.0142$&$-0.0141$\\\bottomrule
\end{tabular}
\caption{Four exact rational certificates for the finite positive range.}\label{tab:positive-certificates}
\end{table}
After substituting the boundary $z=C_0p^2$, the reduced expression for $E_0$ increases with $p$, whereas that for $E_1$ decreases. Thus only the upper and lower endpoint of each numerator block are required. This is distinct from the decrease in $z$ with $p$ fixed used above.

For $p=1,2$, the inequalities $q>C_0p^3$ imply respectively $q\ge9,65$, which lie in the ranges already covered by \cite{RV05}. In every remaining case, \eqref{eq:positive-integer} gives positive linear forms tending to zero along odd $n$. Lemma~\ref{lem:linearforms} completes the proof of Theorem~\ref{thm:positive}.

\section{The negative threshold}\label{sec:negative}
The aim of this section is to establish a common decay estimate for three adjacent negative-argument forms and optimize the resulting condition.
Put $Z=q/p>1$ and $z=-Z$. Let $P_n$ be the Legendre polynomial normalized by $P_n(1)=1$, and define $W_{n-1}(u)=\sum_{k=1}^nP_{k-1}(u)P_{n-k}(u)/k$ and $u=(z+x)/(z-x)=(Z-x)/(Z+x)$.
For $n\ge1$ set
\begin{equation}\label{eq:negative-forms}
g_n(x)=\frac{x^n}{(z-x)^{n+1}}[-P_n(u)\log x-2W_{n-1}(u)],
\qquad F_{n,r}=\int_0^1(1-x)^{2n+r}g_n(x)\dd x,
\end{equation}
where $r=0,1,2$ and $g_n(0)=0$ by continuity.

\subsection{Arithmetic and a common majorant}
For $z>1$, the inner real integral of the Rhin--Viola family is $2(z-x)^{-n-1}Q_n(u)$, and its normalized residue is $(z-x)^{-n-1}P_n(u)$. These are the classical integral formulas for the Legendre functions \cite[(14.7.4), (14.7.6)--(14.7.7), and (14.12.11)]{DLMF}, obtained by $n$ integrations by parts from Rodrigues's formula. Since $2Q_n(u)=P_n(u)\log((u+1)/(u-1))-2W_{n-1}(u)$ and $(u+1)/(u-1)=z/x$, subtracting the $\log z$ residue term gives $F_{n,r}=z^{2n}I_z(2n+r,n,n,n,n)$. In this corrected Legendre representation the logarithm in the integrand is $\log x$, independent of $z$. Set $t=1/z$. The denominator is a power of $1-tx$, whose reciprocal is analytic as a continuous-function-valued map for $|t|<1$. The factor $x^n\log x$ is continuous at zero. Multiplying the polynomial identity by $t^{2n+r}$ clears its apparent pole at $t=0$, and both sides are real analytic throughout $(-1,1)$. Their agreement for $0<t<1$ extends to $-1<t<0$ by the identity theorem. This also matches the negative-argument normalization of \cite[Theorem~2.1]{RV19}, without a separate $\log(-Z)$ term in \eqref{eq:negative-forms}.

The corresponding invariants in \eqref{eq:rv-input} are $H=K=\delta=2n+r$, $\alpha=2n$ and $\beta=0$. Since $r\le2$, a common denominator gives
\begin{equation}\label{eq:negative-integer}
d_{2n+2}^{\,2}p^{2n+2}F_{n,r}\in\Z+\Z\Li(-p/q),\qquad r=0,1,2.
\end{equation}
Here the degree bound in \eqref{eq:rv-input} is what makes the factor $p^{2n+2}$ sufficient.

For $0\le x\le1$, we have $u\in[-1,1]$. The classical bound $|P_k(u)|\le1$ gives $|W_{n-1}(u)|\le H_n$. Therefore, if $M_-(Z)=\max_{0\le x\le1}x(1-x)^2/(Z+x)$, then
\begin{equation}\label{eq:negative-bound}
|F_{n,r}|\le\frac{M_-(Z)^n}{Z}\int_0^1[-\log x+2H_n]\dd x
=\frac{1+2H_n}{Z}M_-(Z)^n\le\frac{1+2n}{Z}M_-(Z)^n.
\end{equation}
The estimate is uniform in the three shifts; the coarser factor $1+2n$ suffices for the exponential argument.

\Needspace{8\baselineskip}
\subsection{Nonvanishing of three adjacent moments}
\begin{lemma}\label{lem:three-moments}
For each $n\ge3$ and $z<-1$, the three numbers $F_{n,0},F_{n,1},F_{n,2}$ do not all vanish.
\end{lemma}
\begin{proof}
Fix $n,z$ and write $\mu_H=\int_0^1(1-x)^Hg_n(x)\dd x$. With $J=x^{-n}g_n$, the Legendre equation and the change of variable $u=(z+x)/(z-x)$ give
$x(z-x)J''+[z-(2n+3)x]J'-(n+1)^2J=0$.
To obtain this equation, substitute the logarithmic representation of $Q_n$ into its ordinary Legendre equation and use the equation for $P_n$. The resulting inhomogeneous equation for $W_{n-1}$ is a polynomial identity, so it extends from $u>1$ to all real $u$. In the expression for $J$, the logarithmic cross terms then cancel. Conjugating by $x^n$ gives
\begin{equation}\label{eq:g-ode}
x^2(z-x)g_n''+[(1-2n)zx-3x^2]g_n'+(n^2z-x)g_n=0.
\end{equation}
For $H\ge2$, integrate \eqref{eq:g-ode} against $(1-x)^H$ and transfer the derivatives to the polynomial coefficients. The resulting recurrence is
\begin{align}\label{eq:moment-recurrence}
0={}&H(H-1)(z-1)\mu_{H-2}
+H[3H-(2H+2n+1)z]\mu_{H-1}\notag\\
&+[z(H+n+1)^2-3H(H+1)-1]\mu_H+(H+1)^2\mu_{H+1}.
\end{align}
For $n\ge3$, the continuous extension of $g_n$ at zero is twice continuously differentiable: its only potentially singular factor is $x^n\log x$. At zero the polynomial coefficients and the required first-order boundary combinations vanish; at one the weight and its first derivative vanish. Thus every boundary term in the two integrations by parts is zero.

If $\mu_{2n}=\mu_{2n+1}=\mu_{2n+2}=0$, the nonzero coefficient $(H+1)^2$ in \eqref{eq:moment-recurrence} propagates the vanishing to all $H\ge2n$. Hence the continuous function $t^{2n}g_n(1-t)$ has zero integral against every polynomial. Uniform polynomial approximation, applied to this function itself, implies that the integral of its square is zero, so it vanishes identically. On the other hand, $P_n(1)=1$ and $W_{n-1}(1)=H_n$ give
$g_n(x)=z^{-n-1}x^n[-\log x-2H_n+O(x|\log x|)]$ as $x\downarrow0$. This function is nonzero for all sufficiently small positive $x$, a contradiction.
\end{proof}

\subsection{Optimization and completion}
The $n$th root of the common factor $d_{2n+2}^2p^{2n+2}$ tends to $e^4p^2$. Thus \eqref{eq:negative-bound}, Lemma~\ref{lem:three-moments} and Lemma~\ref{lem:linearforms} prove irrationality whenever $e^4p^2M_-(Z)<1$.
Put $A=e^4p^2$. This condition is equivalent to $Z>\max_{[0,1]}[Ax(1-x)^2-x]$. Since $A>9$, the maximum is attained at $x_*=(2-\sqrt{1+3/A})/3$ and equals
\begin{equation}\label{eq:maxnegative}
\frac{2A}{27}\left[1+\left(1+\frac3A\right)^{3/2}\right]-\frac23.
\end{equation}
To check global maximality, differentiate the cubic, compare its two critical points with $[0,1]$, and compare the interior value with the endpoint values $0$ and $-1$. Substitution gives \eqref{eq:maxnegative}; it is greater than $1$ for $A\ge e^4$. Consequently the condition itself ensures $Z>1$.

Multiplying \eqref{eq:maxnegative} by $p$ gives \eqref{eq:negative-threshold}. Finally, $(1+t)^{3/2}<1+3t/2+3t^2/8$ for $t>0$: the third derivative of $(1+t)^{3/2}$ is negative, so the quadratic Taylor polynomial at zero is a strict upper bound. This proves the first simpler bound. The correction $-p/3+1/(4e^4p)$ is negative for $p\ge1$, proving the second. Theorem~\ref{thm:negative} follows.

\section{Interpolation for fixed numerators}\label{sec:interpolation}
This section derives the first four rows of Table~\ref{tab:thresholds} from existing certified profiles and the convexity of positive partition sums.
For each $n$, the partition in \cite{Liu} can be written as a summable series $\mathcal Z_n(t)=\sum_\lambda a_{n,\lambda}t^{w_{n,\lambda}}$, with $a_{n,\lambda}\ge0$ and $w_{n,\lambda}\ge0$. H\"older's inequality therefore gives
\begin{equation}\label{eq:holder}
\mathcal Z_n(t_0^{1-\theta}t_1^\theta)
\le\mathcal Z_n(t_0)^{1-\theta}\mathcal Z_n(t_1)^\theta,
\qquad 0\le\theta\le1.
\end{equation}
Together with monotonicity in $t$, this remains an upper bound whenever $t\le t_0^{1-\theta}t_1^\theta$. Endpoint rates $B_0,B_1$ consequently give the rate $(1-\theta)B_0+\theta B_1$ at $t$.

The four rational choices in Table~\ref{tab:interpolation} satisfy the required logarithmic inequality. The endpoint rates are the exact rational bounds extracted from the corresponding certificates in \cite{Liu}. Their source identifiers and exact inequalities are retained in the ancillary files.
\begin{table}[ht]
\centering\small
\begin{tabular}{rrrrrr}\toprule
$t_0$&$t_1$&$\theta$&$t$&usable rate $B$&height cutoff $C$\\\midrule
$3/56$&$2/17$&$23/1000$&$3/55$&$-193/100$&$197/200$\\
$3/62$&$3/56$&$4/25$&$3/61$&$-1057/500$&$451/500$\\
$4/183$&$3/62$&$101/1000$&$4/169$&$-3251/1000$&$99/200$\\
$4/183$&$3/62$&$7/200$&$4/178$&$-3333/1000$&$93/200$\\\bottomrule
\end{tabular}
\caption{Interpolation data. The second and fourth rows are used for negative values of the dilogarithm.}\label{tab:interpolation}
\end{table}
For the indicated sign and target $t=p/Q$, exact logarithmic bounds prove
$B+\log(Q/p)<C<R_\sigma(p,Q)+\log(Q/p)$.
For $q\ge Q\ge3p$, the function $R_\sigma(p,q)+\log(q/p)$ is increasing, so Proposition~\ref{prop:partition} gives the entire denominator range. Its derivative is $(17q-44p)/(12q(q-p))$ for the positive sign and $(17q+44p)/(12q(q+p))$ for the negative sign.

The preceding endpoints are improved to $54,60,168$ by $3/54=1/18$, $-3/60=-1/20$ and $4/168=1/42$. Similarly, the source thresholds $13$ and $17$ for numerators $2$ and $-2$ include $12$ and $16$ after $2/12=1/6$ and $-2/16=-1/8$. The unit-fraction results in the pinned source \cite{Liu}, together with the classical values at $\pm1$, supply the first row. This establishes the first four rows of Table~\ref{tab:thresholds}.

\section{A quadratic comparison measure}\label{sec:quadratic}
The aim of this section is to obtain an analytic rate below $-4.2361$ at $t=5/418$ and transfer it to the last three rows of Table~\ref{tab:thresholds}.
Set $a=12/25$, $L=121/50$, $D=L-a=97/50$, and $w=9/20$. Define $k_1=2(1-a)(1-w)/D^2=1430/9409$ and $k_2=3(1-a)w/D^3=87750/912673$, and let
\begin{equation}\label{eq:density}
\rho(x)=\begin{cases}1,&0\le x\le a,\\
k_1(L-x)+k_2(L-x)^2,&a<x\le L,\\
0,&\text{otherwise}.
\end{cases}
\end{equation}
The tail mass is $1-a$, and its maximum is $k_1D+k_2D^2<1$. Thus $\rho$ is a probability density bounded by $1$.
Define its logarithmic potential and energy by $U(x)=\int\rho(t)\log|x-t|\dd t$ and $E=\int U(x)\rho(x)\dd x$. Put
\begin{align}\label{eq:profile}
H(x)&=3(x+1)\log(x+1)-2x\log x-(x+4)\log(x+4)+1,\notag\\
f(x)&=2U(x)+\tfrac12H(x)-\tfrac12x\log(418/5),\qquad d=-539387/100000,\notag\\
M&=d+\int_0^\infty\pos{f(x)-d}\dd x,
\end{align}
with the usual continuous convention $0\log0=0$.

\subsection{Exact potential and energy certificates}
\begin{lemma}\label{lem:quadratic-numeric}
For the density \eqref{eq:density}, $E>-1.02905$ and $M<-5.16753$. In particular,
\begin{equation}\label{eq:quadratic-rate}
12\log2+4+4(M-E)<-4.2361.
\end{equation}
\end{lemma}
\begin{proof}
The potential, its primitive, and the energy are finite rational linear combinations of powers and logarithms at rational arguments; the formulas and certification method are given in Appendix~\ref{app:certificates}. The positive part in \eqref{eq:profile} is certified by five rational root boxes and a complete sign cover of their complement. On each root box we use the coarse bound $\pos{f-d}<32$. Their total length is $15\cdot10^{-8}$, so their combined contribution is less than $4.8\cdot10^{-6}$. Integrating the primitive on the three certified positive components and adding this error gives the asserted bound for $M$. The energy formula gives the stated lower bound for $E$. All logarithmic bounds are obtained from rational series with explicit remainders. Finally $12\log2+4+4(M-E)<-4.2361$ follows by rational interval arithmetic.
\end{proof}

\subsection{Passage to finite circles}
The finite-circle comparison theorem of \cite{Liu} applies to disjoint circles with positive rational centres. The next argument supplies such a fixed family from \eqref{eq:density}.

\begin{proposition}\label{prop:quadratic-partition}
For all sufficiently large $n$, $\mathcal Z_n(5/418)\le\exp(-4.236n^2)$.
\end{proposition}
\begin{proof}
For a positive integer $N$, take the $N$ midpoint quantiles $x_j$ of $\rho(x)\dd x$. Since $\rho\le1$, consecutive quantiles are separated by at least $1/N$ and $x_1\ge1/(2N)$. Choose rational centres $c_j$ with $|c_j-x_j|<1/(100N)$ and a common radius $r_N=1/(4N)$. The centres are positive and bounded by $16$, the circles are disjoint, and their minimum separation is at least $0.98/N$. Assign equal mass to the circles. Their potential on the real axis is
\begin{equation}\label{eq:discrete-potential}
U_N(x)=\frac1N\sum_{j=1}^N\log\max(r_N,|x-c_j|).
\end{equation}

We use the subsequence $N=m^2$ and put $\delta=1/m$. The kernel $\log\max(\delta,|x-t|)$ is $1/\delta$-Lipschitz in $t$, uniformly for real $x$. Midpoint-quantile quadrature and the rational perturbations therefore give a uniform error at most $(L+0.01)/(N\delta)$. Removing the truncation changes the continuous potential by at most $2\delta$, because $\rho\le1$ and $\int_{-\delta}^{\delta}\log(\delta/|t|)\dd t=2\delta$.

For the discrete potential, choose a centre nearest to $x$. Every other centre, at index distance $j$, is at distance at least $0.49j/N$ from $x$. The inequality $\log_+t\le t$ gives the uniform truncation error bound
\begin{equation}\label{eq:discrete-error}
\frac{\log(4N\delta)}N+\frac{2\delta}{0.49}H_N.
\end{equation}
Since $H_N\le1+\log N$, all three errors tend to zero for $N=m^2$. Thus $U_N\to U$ uniformly on the entire real axis. By the mean-value identity for logarithmic potentials of disjoint circles, their energy is $E_N=N^{-1}\sum_jU_N(c_j)$. Uniform convergence and quantile quadrature for the continuous function $U$ give $E_N\to E$.

The exact inequalities $2/(4-L)-\tfrac12\log(418/5)<0$ and $f(4)<d$ imply that $f$ decreases on $[4,\infty)$. Indeed, $H$ is decreasing there, and for $x\ge y\ge4$ the logarithmic inequality $\log u-\log v\le(u-v)/v$ gives $U(x)-U(y)\le(x-y)/(4-L)$. It follows that $g=\pos{f-d}$ is continuous and supported in $[0,4]$ on the nonnegative axis.

Fix $\varepsilon>0$ small enough to preserve the strict gap between \eqref{eq:quadratic-rate} and $-4.236$. Choose and fix one sufficiently large $N=m^2$ for which both the potential and energy errors are smaller than $\varepsilon$. The corresponding profile satisfies $f_N\le d+g+2\varepsilon$ on the entire nonnegative axis. Since $g$ is continuous and compactly supported, its equally spaced Riemann sums converge to its integral. Each positive integer node occurs at most twice among the $2n$ effective nodes of the determinant construction. Hence their total profile contribution is at most $2n(M+3\varepsilon)$ for all sufficiently large $n$.

Substitute this estimate in the finite-circle comparison inequality of \cite{Liu}. Its remaining errors are $O(n\log n)$ for the fixed family of circles. The choice of $\varepsilon$ and the strict rate gap give the asserted eventual bound. The order of choices is essential: $N$ is fixed before the determinant order $n$ tends to infinity.
\end{proof}

\subsection{Endpoint inequalities and complete denominator ranges}
With $B=-4.236$, rational logarithmic bounds show
$B+\log(418/5)<1901/10000<R_\sigma(p,Q)+\log(Q/p)$
at the six triples
\begin{equation}\label{eq:six-endpoints}
(p,Q,\sigma)=(5,418,1),(5,434,-1),(6,1039,1),(6,1058,-1),(7,2234,1),(7,2256,-1).
\end{equation}
Since $p/Q\le5/418$, Proposition~\ref{prop:partition} applies. The monotonicity used in Section~\ref{sec:interpolation} then gives every denominator beyond each endpoint. The equalities $6/1038=1/173$ and $7/2233=1/319$ add the two preceding positive endpoints. This proves the remaining rows of Table~\ref{tab:thresholds}. The additional rows of Theorem~\ref{thm:finite} follow from the cubic bounds and the exact reductions described in Appendix~\ref{app:certificates}.

\appendix
\section{Exact certificates and reproducibility}\label{app:certificates}
The aim of this appendix is to describe the finite data and exact arithmetic supporting the numerical inequalities in the proofs.
All displayed terminating decimals in inequalities are interpreted as rational numbers. To bound a logarithm, reduce its argument to a fixed compact interval by extracting powers of $2$, and use
$\log x=2\sum_{j=0}^{m-1}v^{2j+1}/(2j+1)+R_m$, where $v=(x-1)/(x+1)$ and $|R_m|\le2|v|^{2m+1}/((2m+1)(1-v^2))$.
This gives explicit rational upper and lower bounds. The generators output exact rational inequalities and complete interval covers for independent checking.

\subsection{Quadratic potential formulas}
For $v\ne0$, put $T(v)=v\log|v|-v$, $J(v)=\tfrac12v^2\log|v|-\tfrac34v^2$, $N(v)=\tfrac16v^3\log|v|-\tfrac{11}{36}v^3$, and $O(v)=\tfrac1{24}v^4\log|v|-\tfrac{25}{288}v^4$, extending all four continuously at zero. These satisfy $J'=T$, $N'=J$ and $O'=N$. With $c=k_1D+k_2D^2-1$, integration of \eqref{eq:density} gives
\begin{align*}
U(x)={}&T(x)+cT(x-a)-(k_1+2k_2D)J(x-a)+2k_2N(x-a)\\
&+k_1J(x-L)-2k_2N(x-L),\\
V(x)={}&J(x)+cJ(x-a)-(k_1+2k_2D)N(x-a)+2k_2O(x-a)\\
&+k_1N(x-L)-2k_2O(x-L),
\end{align*}
where $V'=U$. If $W(v)=\tfrac12v^2\log|v|-\tfrac14v^2$, a primitive of $f$ is
$2V(x)+\tfrac12[3W(x+1)-2W(x)-W(x+4)+x]-\tfrac14x^2\log(418/5)$.
These formulas follow by integrating the piecewise quadratic density; the continuous values at $a$ and $L$ fix the constants.

For the energy, put $A(v)=\tfrac13v^3\log|v|-\tfrac49v^3$ and $P(v)=\tfrac14v^4\log|v|-\tfrac5{16}v^4$, with continuous extensions at zero. Thus $A'(v)=vT(v)$ and $P'(v)=v^2T(v)$. Set
\begin{align*}
C_1={}&L[J(L)-J(a)]-[A(L)-A(a)]-DJ(D)+A(D),\\
C_2={}&L^2[J(L)-J(a)]-2L[A(L)-A(a)]+P(L)-P(a)\\
&-D^2J(D)+2DA(D)-P(D).
\end{align*}
The initial interval, the two cross-interactions and the tail self-interaction give
\begin{align*}
E={}&a^2\log a-\tfrac32a^2+2k_1C_1+2k_2C_2\\
&+(1-a)^2\left[\log D-\frac74-\frac w6-\frac{w^2}{12}\right].
\end{align*}
For example, the normalized linear and quadratic tail densities on $[0,1]$ are $2(1-t)$ and $3(1-t)^2$. Their self-energies are $-7/4$ and $-2$, and their mutual energy is $-11/6$; weighting them by $1-w$ and $w$ gives the last bracket. These values follow by splitting the square into $s<t$ and $s>t$ and integrating powers times logarithms. All arguments in the displayed energy formula are positive rationals, so the stated energy interval follows directly from rational logarithmic bounds.

The five root boxes for $f-d$ are $10^{-8}$ times $[41820077,41820080]$, $[47875826,47875829]$, $[79183885,79183888]$, $[139988153,139988156]$, and $[166842325,166842328]$. 
Away from these boxes the positive components are $[0,41820077/10^8]$, $[47875829/10^8,79183885/10^8]$ and \\$[139988156/10^8,166842325/10^8]$. The complementary sign decisions are certified on 167 rational intervals, including the changes in the explicit formulas at $a$ and $L$. The global cover and the tail argument ensure that no omitted component contributes to $M$.

\subsection{Additional threshold rows}
The aim of this table is to record the remaining uniform tails from the complete threshold list. Each entry follows from Theorems~\ref{thm:positive} and~\ref{thm:negative}, together with exact cancellation at the finitely many preceding denominators. The exact reduction identities are listed in the ancillary data.
\begin{table}[ht]
\centering\small
\begin{tabular}{rrrrrr}\toprule
$p$&$Q_+$&$Q_-$&$p$&$Q_+$&$Q_-$\\\midrule
8&4\,142&4\,138&17&39\,740&39\,734\\
9&5\,897&5\,894&18&47\,172&47\,166\\
10&8\,088&8\,084&19&55\,480&55\,474\\
11&10\,766&10\,763&20&64\,708&64\,702\\
12&13\,976&13\,974&25&126\,385&126\,377\\
13&17\,771&17\,767&50&1\,011\,074&1\,011\,060\\
14&22\,196&22\,190&100&8\,088\,614&8\,088\,582\\
15&27\,300&27\,294&1\,000&8\,088\,614\,820&8\,088\,614\,484\\
16&33\,130&33\,126&&&\\
\bottomrule\end{tabular}
\caption{Additional sufficient thresholds, completing the 24-row list.}\label{tab:additional}
\end{table}

\subsection{Contents of the ancillary package}
The package (available at https://zenodo.org/records/23247332) contains the exact positive floor-interval and logarithmic certificates, the polynomial verification of \eqref{eq:moment-recurrence}, the scalar and global sign certificates for \eqref{eq:density}, and the endpoint-rate data used in Table~\ref{tab:interpolation}. It also contains a pinned copy of the source construction with its original license and attribution. Running \texttt{reference/verify\_all.py} checks the original finite certificates. The versioned source and the explicit parameter tables make the dependence of each numerical claim reproducible.

\Needspace{10\baselineskip}
\section*{Acknowledgments and use of artificial intelligence}
We acknowledge the arithmetic integral constructions of Georges Rhin and Carlo Viola, Aabir Fauzan's determinant approach to $\zeta(5)$, Qian Tang's adaptation to the dilogarithm, and Si-Qi Liu's extension to rational arguments. The formal development underlying our fixed-numerator results comes through Liu's \texttt{Li2Rational}, which derives from Tang's \texttt{Li2NegHalf} (version of 27 September 2026). It also includes code adapted by Liu from Moritz Firsching's \texttt{Zeta5} formalization \cite{Firsching}, notably for Andr\'eief and Gram determinant identities, logarithmic potentials and energies, and arithmetic divisibility and nonvanishing. The original licenses, source headers and detailed attributions are retained in the ancillary package. We also acknowledge the Lean and Mathlib contributors whose work supports these developments.

The development of the arguments and estimates specific to this manuscript was carried out primarily by ChatGPT (Codex), working from the cited mathematical results and supplied source material. Its contributions included developing proof strategies, deriving intermediate lemmas, exploring and selecting parameters, constructing and checking exact certificates, writing and debugging the Lean formalization, and drafting and revising the manuscript. The authors' role was to initiate and direct this process: they supplied the problem and source material, set the objectives, and requested successive mathematical developments, Lean verification and revisions. Thus the AI contribution includes substantive mathematical reasoning as well as computation, formalization and exposition. The authors take responsibility for the manuscript and for accurately reporting the provenance and verification status of its results.

The new estimates and deductions have been checked in Lean 4 with Mathlib and the pinned development of \cite{Liu}. The general theorems retain seven explicit published or classical inputs: the least-common-multiple and finite-prime-window consequences of the prime number theorem; two specialized Rhin--Viola polynomial identities, one including their established small-numerator corollaries; the integral and logarithmic formulas for the Legendre function of the second kind; the unit-interval Legendre bound; and the Legendre differential equation. The ancillary dependency document states these inputs precisely. They are theorem parameters, and their proofs are not re-formalized here. The first seven positive and negative threshold rows, including the $418$ partition estimate, are checked using the pinned source without these external parameters. The complete project compiles, and the public-theorem audit reports only the standard foundational axioms of Lean; the new source contains no admitted proofs or additional axioms.


\begin{thebibliography}{99}
\bibitem{Beukers} F. Beukers, A note on the irrationality of $\zeta(2)$ and $\zeta(3)$, \emph{Bull. London Math. Soc.} 11 (1979), 268--272. 
%\href{https://doi.org/10.1112/blms/11.3.268}{doi:10.1112/blms/11.3.268}.
\bibitem{Davenport} H. Davenport, \emph{Multiplicative Number Theory}, 2nd ed., Graduate Texts in Mathematics 74, Springer, New York, 1980. 
%\href{https://doi.org/10.1007/978-1-4757-5927-3}{doi:10.1007/978-1-4757-5927-3}.
\bibitem{Fauzan} A. Fauzan, $\zeta(5)$ is irrational, preprint, Zenodo, version 1, 17 September 2026. \href{https://doi.org/10.5281/zenodo.22826419}{doi:10.5281/zenodo.22826419}.
\bibitem{Firsching} M. Firsching, Zeta5, software source (2026), commit \href{https://github.com/mo271/Zeta5/tree/f19a1960609f7d38e7b63fd2acb05e6f60a7b741}{\texttt{f19a1960609f}}. \url{https://github.com/mo271/Zeta5}.
\bibitem{Hata} M. Hata, Rational approximations to the dilogarithm, \emph{Trans. Amer. Math. Soc.} 336 (1993), 363--387. 
%\href{https://doi.org/10.1090/S0002-9947-1993-1147401-5}{doi:10.1090/S0002-9947-1993-1147401-5}.
\bibitem{RV05} G. Rhin and C. Viola, The permutation group method for the dilogarithm, \emph{Ann. Sc. Norm. Super. Pisa Cl. Sci. (5)} 4 (2005), 389--437. 
%\href{https://doi.org/10.2422/2036-2145.2005.3.02}{doi:10.2422/2036-2145.2005.3.02}.
\bibitem{RV19} G. Rhin and C. Viola, Linear independence of 1, $\operatorname{Li}_1$ and $\operatorname{Li}_2$, \emph{Moscow J. Comb. Number Theory} 8 (2019), 81--96. %\href{https://doi.org/10.2140/moscow.2019.8.81}{doi:10.2140/moscow.2019.8.81}.
\bibitem{Prellberg} T. Prellberg, Three irrationality results for the dilogarithm: Rhin--Viola and Viola--Zudilin constructions at $-1/4$, $1/5$, and $-1/3$, preprint (2026), \href{https://arxiv.org/abs/2609.06895}{arXiv:2609.06895}.
\bibitem{DLMF} NIST Digital Library of Mathematical Functions, Chapters 14 and 18, in particular Sections 14.2, 14.7, 14.12 and 18.14. \url{https://dlmf.nist.gov/14}.
\bibitem{Tang} Q. Tang, Irrationality of $\operatorname{Li}_2(1/2)$, software source (2026). \url{https://github.com/dtq1997/li2-half-irrationality}.
\bibitem{Liu} S.-Q. Liu, Li2Rational, software source, commit \href{https://github.com/siqiliu-tsinghua/Li2Rational/tree/b2a61c46c157860b57faa163341fda725c629139}{\texttt{b2a61c46c157}}. \url{https://github.com/siqiliu-tsinghua/Li2Rational}.
\end{thebibliography}
\end{document}